\documentclass[12pt,twoside]{amsart}\usepackage{amssymb}
\usepackage{enumerate}

\newtheorem{thmspec}{\relax}

\newtheorem{theorem}{Theorem}[section]
\newtheorem{thm}[theorem]{Theorem}
\newtheorem{lem}[theorem]{Lemma}

\newtheorem{prop}[theorem]{Proposition}
\newtheorem{defi}[theorem]{Definition}

\theoremstyle{definition}

\newtheorem{example}[theorem]{Example}

\theoremstyle{remark}

\numberwithin{equation}{section}
\def \Bbb{\mathbb}

\def\onto{{\kern3pt\to\kern-8pt\to\kern3pt}}

\def\<{\langle}
\def\>{\rangle}
\def\|{{\ |\ }}

\def\onto{\twoheadrightarrow}

\def\-{\underline}

\def\N{\Bbb N}

\def\R{\Bbb R}
\def\C{\Bbb C}
\def\P{\Bbb P}

\def\codim{\operatorname{codim}}
\def\Id{\operatorname{Id}}
\def\d{\operatorname{d}}
\def\ddc{\operatorname{dd^c}}
\def\Crit{\operatorname{Crit}}
\def\I{\mathcal{I}}

\def\<{\langle}
\def\>{\rangle}

\catcode`\@=11
\def\serieslogo@{\relax}
\def\@setcopyright{\relax}
\catcode`\@=12

\title[Algebraic degrees for iterates      of meromorphic self-maps of $\P^k$]
{Algebraic degrees for iterates      of   meromorphic self-maps of $\P^k$}

\begin{document}

\author{Vi{\^e}t-Anh  Nguy\^en}
\address{Vi{\^e}t-Anh  Nguy\^en\\
Mathematics Section\\
The Abdus Salam international centre
 for theoretical physics\\
Strada costiera, 11\\
34014 Trieste, Italy}
\email{vnguyen0@ictp.trieste.it}

\date{}


\begin{abstract}
We first  introduce    the class of quasi-algebraically stable  meromorphic maps of $\P^k.$
This class is strictly larger than that of    algebraically stable  meromorphic self-maps of $\P^k.$ 
  Then we prove that all maps in the new class
enjoy  a recurrent property.  In particular,
the algebraic degrees for iterates of these maps can be computed and    their  first dynamical degrees are
 always  algebraic integers.
   \end{abstract}
\maketitle

\section{ Introduction}

    Let $f: \P^k\longrightarrow \P^k$ be a meromorphic self-map. It can be
    written  $f:= \lbrack F\rbrack:= \lbrack F_0:\ldots :F_k\rbrack$ in homogeneous coordinates
    where  the $F_j$'s  are homogeneous polynomials in the $k+1$ variables $z_0,\ldots,z_k$ of the same degree $d$
    with no nontrivial  common factor.
    The polynomial $F$ will be called {\it a lifting} of   $f$ in $\C^{k+1}.$ The number $\d(f):=d$ will be called  {\it the
    algebraic degree} of  $f.$ Moreover $f$ is said to
    be {\it dominating} if it is generically of maximal rank $k, $ in other words,  its jacobian
    determinant  does not vanish identically (in  any local chart).  The   {\it indeterminacy locus} $\I(f)$ of $f$  is the set of all points of $\P^k$ 
    where $f$ is not  holomorphic, or equivalently the common zero set
    of $k+1$ component polynomials $F_0,\ldots,F_k.$  Observe that $\I(f)$ is a subvariety of codimension at least $2.$
    From now on, we always consider   dominating meromorphic self-maps $f$ of $\P^k$ with $k\geq 2.$
    For such a map $f,$  its {\it first dynamical degree} $\lambda_1(f)$ is given by
    \begin{equation}\label{eq1.1}
    \lambda_1(f):=\lim\limits_{n\to\infty} \d(f^n)^{\frac{1}{n}}.
    \end{equation}
  Computing  $\lambda_1(f)$  and other dynamical degrees associated to $f$  (see, for example, \cite{ds3}
  for a definition of the latter degrees)
  is a fundamental problem in Complex Dynamics. Indeed, this is a necessary step 
  to determine the topological entropy of $f$  (see  \cite{ds3}, \cite{gr}).   

Recall    the following definition (see \cite{fs,fs2,si})
\renewcommand{\thethmspec}{Definition 1}\begin{thmspec}
A  meromorphic self-map $f:\P^k\longrightarrow \P^k$   is said to be  algebraically stable (or  AS for short) if
 there is no  hypersurface of $\P^k$ which would be sent, by some
iterate $f^N,$ to    $\I(f).$  
\end{thmspec}

In other words, $f$ is AS if and only if $\d(f^n)=\d(f)^n,$ $n\in\N.$

 For all AS  maps $f$ with $\d(f) > 1,$ Sibony proves  in \cite{si} that
 the following limit in the sense of current
\begin{equation*}
\lim_{n\to\infty} \frac{(f^n)^{*} \omega}{\d(f^n)}
\end{equation*}
exists,
where  $\omega$ denotes the Fubiny-Study K\"{a}hler form on $\P^k$ so normalized that
 $\int\limits_{\P^k}\omega^k=1$ (see Section 2 below).  This limit is called {\it the Green current} associated to $f.$
The Green current contains many important
dynamical informations  of the corresponding map.  We address the reader to the recent works of 
 Sibony \cite{si},  Dinh--Sibony \cite{ds1,ds2},  Forn{\ae}ss-Sibony  \cite{fs,fs2},  and Favre \cite{fa}
   for    further explanations.

However, the situation becomes much harder in the case of non AS maps.
In general, for a non AS map $f,$  $(f^m)^{\ast}\circ (f^n)^{\ast}\not=(f^{m+n})^{\ast}$  when these operators
act on currents.
  One of the first works 
 in this direction is the article of Bonifant--Forn{\ae}ss \cite{bf}   where    some special non AS
maps are thoroughly studied.  In her thesis   \cite{bo} Bonifant constructs an appropriate Green current for these maps
and then writes down the functional equation. 
  Favre and Jonsson \cite{fj} have studied   the case of polynomial
maps in $\C^2.$ These works show that there is a deep connection between 
the construction of a good invariant current for  a non AS meromorphic self-map $f$
and the sequence $(\d(f^n))_{n=1}^{\infty}$   of algebraic degrees for iterates of $f.$ Moreover, this sequence 
    is, in general,  very complicated.

 The purpose of this
paper is to study   the   sequence $(\d(f^n))_{n=1}^{\infty}$ for all  $f$  in a new class
of  meromorphic self-maps of $\P^k:$   the class of {\it quasi-algebraically stable self-maps}
 (or QAS for short).
This class  contains strictly  that of all   AS  self-maps. 
The QAS self-maps share a {\it recurrent property} with the AS ones. Let us explain
this more explicitly. For an AS self-map $f$   we may define a sequence of polynomial maps $(F_n)_{n=0}^{\infty}:\ \C^{k+1}\longrightarrow
\C^{k+1}$ such that $F_n$ is a lifting of $f^n,$ $n\geq 0,$ in the following recurrent way:
\begin{equation*} 
F_n:=F_1\circ F_{n-1},\qquad n\geq 1,
\end{equation*}
where $F_1,\ F_0$ are arbitrarily fixed  liftings of $f,\ f^0:=\Id$ respectively. 
  Following the same pattern,    the {\it recurrent law} for a QAS self-map $f$ which is not  AS may be stated as follows:
  \begin{equation*}
F_n:=\frac{F_1\circ F_{n-1}}{H_0\circ F_{n-n_0-1}}
\end{equation*}
for all $n>n_0.$ Here $n_0\geq 1$ is an integer and  $H_0$ is a   homogeneous polynomial.
The recurrent phenomenon happens when  the orbits of the hypersurfaces which are sent to $\I(f)$ by some
 iterate $f^N$ (see Definition 1 above) are, in some sense,
not so complicated.  That is the main point of our observation.
 
\smallskip
 
This paper is organized as follows.

\smallskip

We begin Section 2 by
collecting some background  
and introducing some notation. This preparatory  is necessary  for us  to
 state the results afterwards.
 
 Section 3 starts with the definition of quasi-algebraically stable    meromorphic self-maps. Then
  we provide some examples  illustrating this definition. 
  
The last section is devoted to the formulation and the proof of  the main theorem. Some examples
 are also analyzed in the light
of this theorem. Finally, we conclude the paper with some remarks. 
 
\medskip

\indent{\it{\bf Acknowledgment.}} The author is supported by The Alexander von Humboldt
Foundation and The Max-Planck Institut f\"{u}r Mathematik in Bonn (Germany) 
during the preparation of this paper. He wishes to express his gratitude to 
these organisations. He
  would like to thank  N. Sibony
 and T.-C. Dinh for introducing him to Complex Dynamics and for many help. 
  He is also   grateful to    C. Favre whose suggestions helped to  clarify the exposition.
 The  first version  of this article   circulated in a form of a preprint in   2002.

\medskip
\section{Background and notation}
   Let  $f$ be  a dominating   meromorphic self-map of $\P^k,$  $\Gamma(f)$  the graph of $f$ in $\P^k\times\P^k$
   and $\pi_1,$ $\pi_2$ the natural projections of $\Gamma(f)\subset \P^k\times\P^k$ onto its factors.
 Then   $\I(f)$ is exactly the set of points $z\in \P^k$ where $\pi_1$ does not admit a local inverse.
    Let $A$ be a (not necessarily irreducible) analytic subset of $\P^k.$
   We define the following analytic sets
   \begin{equation*}
   f(A):=\overline{f|_{\P^k\setminus \I(f)}(A)}=\overline{\pi_2(\pi_1^{-1}(A\setminus \I(f)))} \qquad\text{and}\qquad f^{-1}(A):=\pi_1(\pi_2^{-1}(A)).
   \end{equation*}
   If $g$ is another  dominating   meromorphic self-map of $\P^k,$ the graph $\Gamma(g\circ f)$
   of the composite map $g\circ f$ is the closure of the set
   \begin{equation*}
   \left\lbrace (z,g(f(z))\in  \P^k\times\P^k:\  z\not\in \I(f)\  \text{and}\  f(z)\not\in \I(g)  \right\rbrace.
   \end{equation*}
  Consequently,  we  can define the sequence  $(f^n)_{n=1}^{\infty}$ of $n$-iterates of $f$  
    by the induction formula  
    \begin{equation*}
    \Gamma(f^n):=\Gamma(f^{n-1}\circ f),\qquad n\geq 2. 
   \end{equation*}
   
   In the sequel, $\codim(A)$ denotes the codimension of $A.$  Moreover, we recall that a {\it hypersurface} is an analytic set of pure codimension $1$ in $\P^k.$
    Let $\Crit(f)$ denote   the critical set of $f$
   (i.e. the hypersurface defined outside $\I(f)$ by  the zero set of the jacobian of $f$ in any local coordinates).
  The following result is very useful.
\begin{lem}\label{irreducible}
Let   $f $  be  as above. Then, for every irreducible analytic set $A\subset\P^k$ such that $A\not\subset\I(f),$
  $f(A)$ is  also an irreducible analytic set.
\end{lem}
\begin{proof} Suppose in order to get a contradiction that $f(A)=B_1\cup B_2,$ where $B_1,$ $B_2$ are
analytic sets in $\P^k,$ distinct from $f(A).$  It follows that $(A\cap f^{-1}(B_1))\cup(A\cap f^{-1}(B_2))\subset A$
and the two  sets $A\cap f^{-1}(B_1),$  $A\cap f^{-1}(B_2)$ are distinct analytic  components of $A.$ We therefore  get the desired contradiction. 
 This finishes the
proof.
\end{proof}

We denote by $\mathcal{C}^{+}_1(\P^k)$ the set of
 positive closed currents  of bidegree  $(1,1)$ on $\P^k.$ A current $T\in  \mathcal{C}^{+}_1(\P^k) $ can be written locally
  as $T=\ddc u$ for some plurisubharmonic
 function $u$ (which is called a {\it local potential} of $T$).  The mass of $T$  is defined by 
  $\Vert T\Vert:= \int\limits_{\P^k} T\wedge \omega^{k-1}.$  Fix a point $a\in\P^k$
 and local coordinates sending $a$ to the origin in $\C^k.$  Choose a local plurisubharmonic  potential $u$ for $T$ defined around $0$ in these coordinates.
 We can define the {\it Lelong number} of $u$ at $0$  as follows
 \begin{equation*}
 \nu(u,0):=\max\left\lbrace c\geq 0:\  u(z)\leq c\log\vert z\vert+\mathcal{O}(1)\right\rbrace
 \end{equation*} 
 which is a finite nonnegative real number. We then set  $\nu(T,a):=
 \nu(u,0),$ which does not depend on any choice we made.
 
 For a current $T\in  \mathcal{C}^{+}_1(\P^k) ,$ we use local potentials to define the induced pull-back
 $f^{\ast}T\in  \mathcal{C}^{+}_1(\P^k).$ More precisely, for any $z\in \P^k\setminus\I(f),$  $T$ has a local potential  $u$
  in a neighborhood of $f(z),$ and we define  $\widetilde{f^{\ast}T}:=\ddc (u\circ f)$ in a neighborhood of $z.$
   This yields a well-defined,
  positive closed $(1,1)$-current on the set $\P^k\setminus \I(f).$ Since $\codim(\I(f))>1,$  we can  extend 
   $\widetilde{f^{\ast}T}$ to a  current  $f^{\ast}T\in  \mathcal{C}^{+}_1(\P^k) $ by assigning 
  zero mass on the set $\I(f)$ to the coefficients measures of $\widetilde{f^{\ast}T}.$

 Any hypersurface $\mathcal{H}$ of $\P^k$ defines a current of integration  $[\mathcal{H}]\in  \mathcal{C}^{+}_1(\P^k),$
 and $\Vert [\mathcal{H}] \Vert =\deg(H),$ where $H:\ \C^{k+1}\longrightarrow\C$ is any homogeneous polynomial defining  $\mathcal{H}.$
 Finally, for any current $T\in  \mathcal{C}^{+}_1(\P^k), $  it holds  that
 \begin{equation}\label{eq2.1}
 \Vert f^{\ast} T\Vert=\d(f)\cdot\Vert T\Vert.
 \end{equation}
 For further information on this matter, the reader is invited to consult the works \cite{fs} and \cite{si}.

\section{Quasi-algebraically stable meromorphic self-maps}
In \cite{fs2} Forn{\ae}ss and Sibony establish the following definition.
\begin{defi}\label{degreelowering} A  hypersurface $\mathcal{H}\subset \P^k$
is said to be a {\it degree lowering  hypersurface} of $f$ if, for some (smallest) $n\geq
1,$   $f^n(\mathcal{H})\subset \I(f).$ The integer $n$ is then called the height of $\mathcal{H}.$
\end{defi}

The following proposition gives us the structure of a non AS self-map.
 \begin{prop}\label{structure}
    Let   $f $ be a   meromorphic self-map of $\P^k.$   Then there are exactly an
    integer $M\geq 0,$  $M$ degree lowering  hypersurfaces $\mathcal{H}_j$   with height
    $n_j,$ $j=1,\ldots,M,$   satisfying the following properties: 
     \\
    (i)  all the numbers $n_j$ $,j=1,\ldots,M,$ are pairwise different;\\
    (ii)  $\codim\left(f^m(\mathcal{H}_j)\right) >1$ for  
    $ m=1,\ldots,n_j,$ and $j=1,\ldots,M;$\\
    (iii) for any   degree lowering  irreducible hypersurface $\mathcal{H}$ of $f,$ there are
     integers $n\geq 0$ and  $1\leq j\leq M$ such that $f^n(\mathcal{H})$ is a  hypersurface and 
    $f^n(\mathcal{H})\subset  \mathcal{H}_j$
 \end{prop}
 In particular,    $f$ is AS if and only if  $M=0.$
   \begin{proof}
     First, we give the construction of  $M$ and $\mathcal{H}_j,$ $n_j,$  $j=1\ldots,M.$   To this end observe that
    every hypersurface $\mathcal{H}$ satisfying $\codim\left(f (\mathcal{H} )\right) >1$ should be contained in $\Crit(f).$ 
   Therefore, one takes the family $\mathcal{F}$ of all degree lowering  irreducible components $\mathcal{H}$  of     $\Crit(f)$  such that $\codim \left(f^m(\mathcal{H})\right) >1$
   for $1\leq m\leq n,$  where  $n$ is the  height of $\mathcal{H}.$
   Let $1\leq n_1<\cdots<n_M$ be all the heights of  elements in $\mathcal{F}$   (it is easy to see that $M$ is  finite). 
   Let $\mathcal{H}_j$ be the (finite) union of all  elements in $\mathcal{F}$  with the same height $n_j.$
   Then  properties
   (i) and (ii) are satisfied.  
   
   To prove (iii) let  $\mathcal{H}$ be a degree lowering  irreducible hypersurface  with height $h.$  
    In virtue of Lemma \ref{irreducible},
   let $n$  be the greatest integer  such that  $0\leq n<h$  and $f^n(\mathcal{H})$ is a hypersurface.
    The choice of $n$ implies that  $\codim\left(f^m(\mathcal{H})\right) >1,$  $m=n+1,\ldots,h.$
    Consequently, in virtue of the above construction, we deduce that $f^n(\mathcal{H})\subset  \mathcal{H}_j$ for some $1\leq j\leq M.$
    This proves (iii).
    
    Since  the uniqueness of  $M$ and $\mathcal{H}_j,$ $n_j,$  $j=1,\ldots,M,$     is almost obvious, it
     is therefore left to the  reader. This completes the proof.
\end{proof}
\begin{defi} \label{primitive} Under the hypothesis and  the notation of  Proposition \ref{structure}, for every $ j=1,\ldots,M,$
 $\mathcal{H}_j$ is called the primitive   degree lowering  hypersurface of
 $f$  with the height $n_j.$ 
\end{defi}

We are now able to  define the class of quasi-algebraically stable self-maps.
\begin{defi} \label{QAS} A meromorphic self-map $f $ of $\P^k$ is said to be   quasi-algebraically stable 
(or QAS for short) if either it is AS or it satisfies the following properties:
 \begin{itemize}
\item[$(i)$]   there is only one primitive degree lowering
hypersurface (let  $\mathcal{H}_0$ be this hypersurface and let  $n_0$ be its height);
\item[$(ii)$]   for every   irreducible component $\mathcal{H}$
of $\mathcal{H}_0$ and every $m=1,\ldots,n_0,$  $f^{m}(\mathcal{H})\not\subset\mathcal{H}_0;$
  \item[$(iii)$] for every   irreducible component $\mathcal{H}$
of $\mathcal{H}_0,$   one of the following two conditions holds\\ 
    $(iii)_1$  $f^{m}(\mathcal{H})\not\subset \I(f)$ for all $m\geq n_0+1,$\\
     $(iii)_2$    there is an $m_0\geq n_0$ such that
 $f^{m_0+1}(\mathcal{H})$ is a hypersurface and $f^{m}(\mathcal{H})\not\subset \I(f)$ for all $m$ verifying $n_0+1\leq m\leq m_0.$
\end{itemize}
\end{defi}
It is worthy to remark that Proposition \ref{structure} allows us to check if a map  is QAS. 

We conclude  this section by studying some examples.
\begin{example}  \label{ex1}
Consider the following meromorphic self-map of $\P^2:$
\begin{equation}\label{eqex1}
 f\left(\left\lbrack z:w:t\right\rbrack\right ):=
 \left \lbrack    2tz- \left(z^2+w^2\right )  :
2tw-\left (  z^2+w^2 \right ) :  2t^2-\left (  z^2+w^2 \right )
\right\rbrack .
\end{equation}
It can be checked that
$
\I(f)=\left\lbrace \left\lbrack
1:1:1\right\rbrack, \left\lbrack 1:i:0\right\rbrack,\left\lbrack
1:-i:0\right\rbrack
\right\rbrace,
$
and   $\Crit(f)=\lbrace t(2t^2+w^2+z^2-2zt-2wt)=0
\rbrace.$
Moreover we have
\begin{eqnarray*}
f\left (\left\lbrace   t=0   \right\rbrace \right )&=&\lbrack
1:1:1\rbrack\in \I(f),\\
 f\left (\left\lbrace   2t^2+w^2+z^2-2zt-2wt =0   \right\rbrace \right )&=&\left\lbrace
t-z-w=0\right\rbrace.
\end{eqnarray*}
 Therefore, $\left\lbrace   t=0   \right\rbrace
$ is the unique primitive  degree lowering hypersurface and its height is  $1.$  
Since $\left\lbrack
1:1:1\right\rbrack\not\in \left\lbrace   t=0   \right\rbrace$  and $f^2(\{t=0\})$ is a hypersurface,     $f$   is QAS.  
 \end{example}

\begin{example}   For all integers $d\geq 2$ and $m\geq 1,$  the following map   is given by Bonifant-Forn{\ae}ss in 
\cite{bf}
\begin{equation*} 
 f\left(\left\lbrack z:w:t\right\rbrack\right ):=
 \left \lbrack    zt^{d-1}  :
\left (  wt^{d-1}+z^d \right )\cos{\frac{\pi}{m}}-t^d \sin{\frac{\pi}{m}}
     : \left (  wt^{d-1}+z^d \right )\sin{\frac{\pi}{m}}+t^d \cos{\frac{\pi}{m}}
\right\rbrack .
\end{equation*}
 It can be checked that $\I(f)=\left\lbrack
0:1:0\right\rbrack,$ $\Crit(f)=\lbrace t=0\rbrace,$ and
$\lbrace t=0\rbrace$ is the only primitive degree lowering hypersurface of $f.$ Moreover, its height  is 
$m.$   Since  $ f^n(\lbrace t=0\rbrace)=\left\lbrack 0:\cos{\frac{n\pi}{m}}:
\sin{\frac{n\pi}{m}}\right\rbrack\in \{t=0\}$  for $n=1,\ldots,m,$ $f$ is not a QAS according to
Definition \ref{QAS} (ii). However, $f$ satisfies conditions (i) and $\text{(iii)}_{\text{2}} $ of this definition.  
\end{example}   
 \begin{example}  \label{ex3}
Consider the following meromorphic self-map of $\P^2:$
\begin{multline}\label{eqex3}
 f\left(\left\lbrack z:w:t\right\rbrack\right ):=
 \Big \lbrack   (z+w+t)^2 (z^3+w^3+t^3)z^2-27z^3w^4   :\\
   (z+w+t)^2 (z^3+w^3+t^3)w^2-27z^3w^4  :   (z+w+t)^2 (z^3+w^3+t^3)t^2-27z^3w^4 
\Big\rbrack .
\end{multline}
It can be checked that
\begin{equation*}
\I(f)=  \left\lbrack
1:1:1\right\rbrack  \cup \left\lbrace [z:w:t],\  z+w+t=zw=0
\right\rbrace\cup \left\lbrace [z:w:t],\  z^3+w^3+t^3=zw=0
\right\rbrace,
\end{equation*}
and    
\begin{equation*}
f\left (\left\lbrace  z+w+t=0   \right\rbrace \right )=f\left (\left\lbrace  z^3+w^3+t^3=0   \right\rbrace \right )=\lbrack
1:1:1\rbrack\in \I(f).
  \end{equation*}
  
  Let $G$ be an irreducible homogeneous polynomial in $\C^3$ such that
  $f(\mathcal{G})$ is a point $[a:b:c]\in\P^2,$ where $\mathcal{G}$ is the hypersurface $\{G =0\}$ in $\P^2.$
  We deduce from (\ref{eqex3}) and the equality  $f(\mathcal{G})=[a:b:c]$  that
  $G$ divides both polynomials 
 $  (z+w+t)^2 (z^3+w^3+t^3)(bz^2-aw^2)-27z^3w^4(b-a)$ and
  $  (z+w+t)^2 (z^3+w^3+t^3)(cz^2-at^2)-27z^3w^4(c-a).$ Hence,
  $G$ divides  the polynomial
  \begin{equation*}
   (z+w+t)^2 (z^3+w^3+t^3)\Big( (c-a)(bz^2-aw^2)-(b-a)(cz^2-at^2)\Big).
  \end{equation*}
   Now it is not difficult to see that there is no   hypersurface $\mathcal{G}$  which 
 is not contained in  $\mathcal{H}_0:=\left\lbrace (z+w+t)^2( z^3+w^3+t^3)=0   \right\rbrace$ and which satisfies
  $\codim(f(\mathcal{G}))>1.$
 In other words, $\mathcal{H}_0$ is  the unique primitive  degree lowering hypersurface and its height is $1.$
Since $\left\lbrack
1:1:1\right\rbrack\not\in  \mathcal{H}_0,$  and  $f^2(\lbrace z+w+t=0 \rbrace),$ $f^2(\lbrace z^3+w^3+t^3=0 \rbrace)$ are hypersurfaces, it follows that     $f$   is QAS.  
 \end{example}

 \section{The main result} 
Now we are ready to formulate the main result of this article.
\renewcommand{\thethmspec}{The Main Theorem}\begin{thmspec}
    \label{mainthm}
    Let  $f$  be a QAS meromorphic self-map of $\P^k$ which is not AS.
    Let $\mathcal{H}_0$ be its  unique primitive degree lowering hypersurface and let $n_0$ be its   height. 
      We define a sequence $\lbrace F_n\ :n\geq 1\rbrace$ of maps  $\C^{k+1} \longrightarrow\C^{k+1}$ as follows~:

      $F_1 ,\ldots,F_{n_0},F_{n_0+1} $ are arbitrarily fixed
liftings of $f^1 (\equiv f),\ldots,f^{n_0},
f^{n_0+1}$ respectively.
Let $H_0$ be the unique  homogeneous polynomial  which verifies the equality
\begin{equation}\label{eq4.1}
F_1\circ F_{n_0}= H_0\cdot F_{n_0+1}.
\end{equation}
Next we   define $F_n$ for all $ n>
n_0+1$ as follows~:
\begin{equation}\label{eq4.2}
F_n:=\frac{F_1\circ F_{n-1}}{H_0\circ F_{n-n_0-1}}.
\end{equation}
 
Then $\mathcal{H}_0=\left\lbrace H_0(z)=0\right\rbrace$, and for any $n\geq 0,$
  $F_n$ is   a lifting of   $f^n.$   Moreover,
  for any  current $T\in\mathcal{C}^{+}_1(\P^k),$   
 \begin{equation}\label{eq4.3}
  (f^n)^{*}T=
\begin{cases}
(f^{n-1})^{*}(f^{*}T), 
  &n=1,\ldots,n_0  ,\\
 (f^{n-1})^{*}(f^{*}T)-\Vert T\Vert\cdot (f^{n-n_0-1})^{*}[\mathcal{H}_0], 
&  n>n_0.
\end{cases}
 \end{equation}
 \end{thmspec}
 
 Prior to the proof of the theorem we need a preparatory result.

\begin{lem}\label{recurrent}  We keep the above hypothesis and notation.
Let $m\in\N,$ $m\geq 1.$ Then, for any current $T:=[l(z)=0],$ where $\{l(z)=0\}$ is  a generic 
(in the  Zariski  sense)
complex hyperplane in $\P^k,$ 
the supports of   $ (f^p)^{\ast}[\mathcal{H}_0]$ and  $ (f^m)^{\ast} (f^{n_0+1})^{\ast}T$ do not
contain any component of  the hypersurface $ (f^m)^{-1}(\mathcal{H}_0) $
for all $p$ verifying $ \max\{0,m-n_0+1\}\leq p\leq m-1.$
\end{lem}
 
  \smallskip
  
  \noindent{\it Proof of Lemma \ref{recurrent}.}
  To prove  the assertion for  $ (f^p)^{\ast}[\mathcal{H}_0]$, fix an arbitrary irreducible component $\mathcal{H}$ of the hypersurface $ (f^m)^{-1}(\mathcal{H}_0), $ and
 an integer $p:$  $\max\{0,m-n_0+1\}\leq p\leq m-1.$ Suppose in order to reach a contradiction that
 $\mathcal{H}\subset (f^p)^{-1}(\mathcal{H}_0).$
 Putting
  $\mathcal{G}:=f^p(\mathcal{H}),$
  the latter inclusion implies that 
  \begin{equation}\label{eqrecurrent3}
  \mathcal{G}\subset\mathcal{H}_0.
  \end{equation}
  
  In virtue of Lemma \ref{irreducible}, there are two cases to consider.
  \\
  {\bf Case 1:} $\mathcal{G}$ is an irreducible hypersurface.
  
  In this case  it follows from the inclusion  $\mathcal{H} \subset (f^m)^{-1}(\mathcal{H}_0) $
  and the equality  $\mathcal{G}=f^p(\mathcal{H})$  that $ f^{m-p}(\mathcal{G})=f^m(\mathcal{H})\subset \mathcal{H}_0.$  
  Recall that $m-p\leq n_0, $  $\mathcal{G}$ is an irreducible component of $ \mathcal{H}_0$   (by (\ref{eqrecurrent3}))
   and $f$ is QAS. Consequently, if follows from Definition \ref{QAS} (ii) that 
     $ f^{m-p}(\mathcal{G})\not\subset \mathcal{H}_0,$ 
  which contradicts the inclusion $ f^{m-p}(\mathcal{G})\subset \mathcal{H}_0.$
  Hence, this case cannot happen. 
  \\
  {\bf Case 2:} $\mathcal{G}$ is an irreducible analytic set of codimension strictly greater than $1.$
  
  Let $q$ be the greatest integer such that $0\leq q<p$ and $f^q(\mathcal{H})$ is a hypersurface. 
  Consider three subcases. 
  
    {\bf Case 2a:} there is a smallest integer $r$ such that  $q<r <p$ and  $f^r(\mathcal{H})\subset\I(f).$ 
    
    In virtue of the choice of $q,$ $r,$  of Definition \ref{QAS} and Lemma \ref{irreducible},
     we see that $f^q(\mathcal{H})$
     is an irreducible component of $ \mathcal{H}_0.$ Next, we will analyze the
     orbit of $f^q(\mathcal{H})$ under iterations of $f.$
     
    Using the choice of $q,$ we see that  the following analytic sets
     $f^{q+1}( \mathcal{H}),\ldots, f^{p-1}(\mathcal{H})$ are of codimension strictly greater than $1.$
     Recall from  the assumption of Case 2 and Case 2a that $f^p(\mathcal{H})(=\mathcal{G})$ is also
      of codimension strictly greater than $1$ and   $f^r(\mathcal{H})\subset\I(f).$
    Consequently, using  Definition \ref{QAS} (iii) we see that none of the following analytic sets
    $f^{q+1}( \mathcal{H}),\ldots, f^{r-1}(\mathcal{H})$ and  $f^{r+1}( \mathcal{H}),\ldots, f^p(\mathcal{H})$
    is contained in $\I(f).$ Hence,
    \begin{equation*}
    f^{p+1}(\mathcal{H})=f(f^p(\mathcal{H})) =f(\mathcal{G})\subset f(\mathcal{H}_0),
    \end{equation*}
    where the inclusion follows from  (\ref{eqrecurrent3}). 
      Since $f$ is QAS, we deduce from the latter  inclusion (i.e. $ f^{p+1}(\mathcal{H})\subset f(\mathcal{H}_0)$)
      and Definition \ref{QAS} (ii)      that there is a smallest integer $s$
    such that  $p<s\leq p+n_0,$ and   $\codim(f^{p+1}( \mathcal{H})),\ldots, \codim(f^{s}(\mathcal{H}))>1, $  
    and none   of the following   sets $f^{t}( \mathcal{H})$ $ (p+1\leq t\leq s-1) $
    is contained in $\I(f),$  but $ f^{s}(\mathcal{H})\subset \I(f).$
    
    In summary, we have shown that  for some $q<r<p<s,$
     \begin{itemize}
    \item[$\bullet$] $f^q(\mathcal{H})$
     is an irreducible component of $ \mathcal{H}_0;$
     \item[$\bullet$]  all the analytic sets 
     $f^{t}( \mathcal{H})$  ($q+1\leq t\leq r-1$) are of  codimension strictly greater than $1,$ and
     none   of them  
    is contained in $\I(f);$ 
      \item[$\bullet$]
     $f^r(\mathcal{H})\subset\I(f);$
     \item[$\bullet$]
      all the analytic sets 
     $f^{t}( \mathcal{H})$    ($ r+1\leq t\leq s-1$) are of  codimension strictly greater than $1,$ and
     none   of them  
    is contained in $\I(f);$  
    \item[$\bullet$]
     $ f^{s}(\mathcal{H})\subset \I(f).$
    \end{itemize}
     Since   $f^{t}( \mathcal{H})=f^{t-q}(f^q(\mathcal{H}))$
     for $q+1\leq t\leq s,$   Definition \ref{QAS} (iii)  says that $f^q(\mathcal{H})$
     cannot satisfy all the above $\bullet$.
     Hence, Case 2a cannot happen.
       
     {\bf Case 2b:}   $f^r(\mathcal{H})\not\subset\I(f),$  $r=q+1,\ldots,p-1,$
     but   $f^p(\mathcal{H}) \subset\I(f).$  
     
     Under this assumption  and using the choice of $q$  and  Definition \ref{QAS}, we see that
      $f^q(\mathcal{H})$ is an irreducible component of $ \mathcal{H}_0$
     and $p=q+n_0.$ Therefore, it follows from  (\ref{eqrecurrent3}) that
      $f^{n_0}( f^q(\mathcal{H}))=f^p(\mathcal{H})\subset \mathcal{H}_0.$
     On the other hand, Definition \ref{QAS} (ii) says that  $f^{n_0}( f^q(\mathcal{H}))\not\subset \mathcal{H}_0.$
     Hence, Case 2b is impossible.

     {\bf Case 2c:}   $f^r(\mathcal{H})\not\subset\I(f),$  $r=q+1,\ldots,p.$

     First we claim that  there is a smallest integer $s$ such that $s>q$ and
     $f^s(\mathcal{H})\subset  \I(f).$
Otherwise, using the fact that $f^q(\mathcal{H})$ is a hypersurface  (by the choice of $q$),
 one would deduce that $f^{n+t}(\mathcal{H})=f^{n}(           f^t(\mathcal{H}))$ for all $t\geq q$  and $n\geq 0.$ 
In particular, using  (\ref{eqrecurrent3}) and the identity $\mathcal{G}=f^p(\mathcal{H}) $ we would have
\begin{equation*}
 f^{p+n_0}(\mathcal{H})= f^{n_0}(\mathcal{G}) \subset f^{n_0}(\mathcal{H}_0)\subset \I(f),
\end{equation*}
which  is a contradiction. Hence, the above claim has been proved.

     Using this claim,   the assumption of Case 2c and the choice of $q,$ we see that
      $p<s ,$ 
and all the  analytic sets  $ f^{t}( \mathcal{H}) \ (q+1\leq t\leq s-1)$  are of 
 codimension strictly greater than $1,$ 
      and none   of them
    is contained in $\I(f),$  but $ f^{s}(\mathcal{H})\subset \I(f).$   This implies that
     $f^q(\mathcal{H})$ is an irreducible component of $ \mathcal{H}_0,$
     and $s=q+n_0.$   Since  $ f^{p-q}(f^q(\mathcal{H}))=
                f^p(\mathcal{H})\subset \mathcal{H}_0$  (by (\ref{eqrecurrent3}))
     and $0<p-q<s- q=n_0,$
      we obtain a contradiction with Definition \ref{QAS} (ii).
        Hence, Case 2c cannot happen.
   
  Hence,    the proof of  the assertion for  $ (f^p)^{\ast}[\mathcal{H}_0]$ is complete.
  
  To prove   the assertion for  $ (f^m)^{\ast} (f^{n_0+1})^{\ast}T,$ fix an arbitrary irreducible component $\mathcal{H}$ of the hypersurface $ (f^m)^{-1}(\mathcal{H}_0), $ and
  a current $T:=[l(z)=0],$ where $\{l(z)=0\}$ is  a generic complex hyperplane of $\P^k.$ 
   Suppose in order to reach a contradiction that
 $\nu\Big( (f^m)^{\ast} (f^{n_0+1})^{\ast}T,  z\Big)>0$ for any generic point $z\in\mathcal{H}.$ Putting $\mathcal{G}:=f^m(\mathcal{H}),$
  the latter inequality and the choice of $\mathcal{H}$  imply that 
  \begin{equation}\label{eqrecurrent4}
  \mathcal{G}\subset\mathcal{H}_0\cap \I(f^{n_0+1}).
  \end{equation}
  
  Let $q$ be the greatest integer such that $0\leq q\leq m$ and $f^q(\mathcal{H})$ is a hypersurface. 
 Clearly,  $q<m$  because of  (\ref{eqrecurrent4}):  $\codim(\I(f^{n_0+1}))>1.$  
 In virtue of the choice of $q,$  and of Definition \ref{QAS}, and of  the inclusion 
 $\mathcal{G}=f^m(\mathcal{H})\subset  \I(f^{n_0+1})$ (see (\ref{eqrecurrent4})), we conclude  that $f^q(\mathcal{H})$ is an irreducible component of $ \mathcal{H}_0.$
  Since  $f^{m-q}(f^q(\mathcal{H}))=f^m(\mathcal{H})\subset \mathcal{H}_0$  (see  (\ref{eqrecurrent4})),
      we obtain a contradiction with Definition \ref{QAS} (ii)--(iii).

  Hence, the proof of   the assertion for  $ (f^m)^{\ast} (f^{n_0+1})^{\ast}T$ is finished. This completes the proof of the lemma.
  \hfill $\square$
  
  \smallskip
     
 Now we arrive at
 \subsection*{Proof of the Main Theorem.}
  The assertion $\mathcal{H}_0=\left\lbrace H_0(z)=0\right\rbrace$ follows immediately from
   (\ref{eq4.1}) and the hypothesis on $\mathcal{H}_0$ and $n_0.$ Moreover, the hypothesis of the theorem implies
   that  $F_n$ is a lifting of $f^n$ and (\ref{eq4.3}) is valid for $n=1,\ldots,n_0.$
   
   We will prove (\ref{eq4.3}) and the fact that  $F_n$ is a lifting of $f^n$ by induction on $n\geq n_0+1.$ For $n=n_0+1,$  these  assertions are  immediate consequences of
  (\ref{eq4.1})--(\ref{eq4.2}). Suppose  them true for   $ n-1,$
we like to show them for $n.$

To this end  let $G$  be the homogeneous polynomial given by 
\begin{equation}\label{eq4.4}
G\cdot F_n:= F_1\circ F_{n-1} ,
\end{equation}
and let $\mathcal{G}$ be the hypersurface $\left\lbrace G(z)=0\right\rbrace.$  We may rewrite (\ref{eq4.4})  as
\begin{equation}\label{eq4.5}
 (f^{n-1})^{\ast}(f^{\ast}T)= (f^n)^{\ast}T +[\mathcal{G}],
\end{equation}
for any current $T\in\mathcal{C}^{+}_1(\P^k)$ of mass $1.$  In virtue of (\ref{eq4.2}) and (\ref{eq4.5}), we only need to show that
\begin{equation}\label{eq4.6}
[\mathcal{G}]=(f^{n-n_0-1})^{*}[\mathcal{H}_0].
\end{equation}
One breaks the proof of this identity into two steps.

\smallskip

\noindent {\bf Step I:}  Proof of the inclusion $ \mathcal{G}\subset (f^{n-n_0-1})^{-1}(\mathcal{H}_0).$

Consider an arbitrary  irreducible component $\mathcal{F}$ of $\mathcal{G}.$
Then we deduce from (\ref{eq4.5}) that  
\begin{equation*}
\nu\Big((f^{n-1})^{\ast}(f^{\ast}T),z\Big)>0
\end{equation*}
for   any current $T:=[l(z)=0],$ where $\{l(z)=0\}$ is  a generic complex hyperplane in $\P^k,$ 
and for a generic point $z\in\mathcal{F}.$ Since for any point  $z$ outside $\I(f),$
we can choose $T$ so that
  $ f^{\ast}T $  vanishes in a neighborhood of $z,$ 
it follows from the latter inequality that $(f^{n-1})(\mathcal{F})\subset \I(f).$

Now   let $ m$   be the greatest  integer
 such that $0\leq m<n-1$ and $f^{m}(\mathcal{F})$ is a hypersurface. Put $\mathcal{H}:=  f^{m}(\mathcal{F}).$
   Therefore,
$f^{m+1}(\mathcal{F}),\ldots,f^{n-1}(\mathcal{F})$ are  analytic sets of codimension strictly greater than $1.$
Since we have shown that $f^{n-1}(\mathcal{F})\subset \I(f),$ there is a smallest integer $p$ such that
$m+1\leq p\leq n-1$  and  $f^p(\mathcal{F})\subset \I(f).$  Using the hypothesis that $f$ is QAS,   one concludes that  
$\mathcal{H}:=f^{m}(\mathcal{F})$ is an irreducible component of  $\mathcal{H}_0.$  Moreover, one has  $p=n_0+m.$

Recall from the previous paragraph that all the analytic sets
$f^t(\mathcal{H})$   ($1\leq t\leq n-1-m$)  are of  codimension strictly greater than $1.$
In addition, we have  that   
  $f^{p-m}(\mathcal{H})=f^p(\mathcal{F})\subset \I(f),$  and  $f^{n-1-m} ( \mathcal{H})=f^{n-1}(\mathcal{F}) 
\subset \I(f)$ with $n-1-m\geq p- m=n_0.$  Invoking Definition \ref{QAS} (iii) and the choice of $p,$ it follows that
    $p-m=n-1-m.$ This, combined with the equality $p=n_0+m,$ implies that $m=n-n_0-1.$  In summary, we have shown that
$f^{n-n_0-1}(\mathcal{F})=\mathcal{H}\subset \mathcal{H}_0.$ Since  $\mathcal{F}$ is an arbitrary component
of $\mathcal{G},$  we deduce that $\mathcal{G}\subset (f^{n-n_0-1})^{-1}(\mathcal{H}_0).$  This completes Step I.

\smallskip

\noindent {\bf Step II:} Proof of  identity (\ref{eq4.6}).\footnote{ The author thanks C. Favre for
 suggesting him the use of currents in this step. This helps to clarify the author's proof.}

  In what follows $T$ is a  current  in $\mathcal{C}^{+}_1(\P^k)$  of mass 1,  and $d:=\d(f).$  Moreover, we make the following convention
   $(f^{ m})^{\ast}[\mathcal{H}_0]:=0$ for all $m<0.$  Next,
  we apply the hypothesis of induction (i.e. identity (\ref{eq4.3})) for $n-1,\ldots, n-n_0$ repeatedly by taking into account the identity
   $\Vert (f^{m})^{\ast} T\Vert=d^m,$ $m=1,\ldots,n_0,$  (see (\ref{eq2.1})).  Consequently,  one gets 
 \begin{equation}\label{eq4.7}
\begin{split}
  (f^{n-1})^{\ast}(f^{\ast})T&= (f^{n-2})^{\ast}(f^{\ast}f^{\ast})T  -d(f^{n-n_0-2})^{\ast}[\mathcal{H}_0]\\
  &=\cdots\\
   &=(f^{n-n_0 })^{\ast} (f^{ n_0})^{\ast} T-d^{ n_0-1} (f^{n-2n_0})^{\ast}[\mathcal{H}_0]-\cdots-d(f^{n-n_0-2})^{\ast}[\mathcal{H}_0]\\
   &=(f^{n-n_0-1})^{\ast}f^{\ast}(f^{n_0})^{\ast}T-       d^{ n_0} (f^{n-2n_0+1})^{\ast}[\mathcal{H}_0]-\cdots-d(f^{n-n_0-2})^{\ast}[\mathcal{H}_0]\\
   &=(f^{n-n_0-1})^{\ast} (f^{n_0+1})^{\ast}T + (f^{n-n_0-1})^{\ast}[\mathcal{H}_0]\\
   &\quad-d^{ n_0} (f^{n-2n_0+1})^{\ast}[\mathcal{H}_0] -\cdots-d(f^{n-n_0-2})^{\ast}[\mathcal{H}_0].
\end{split}\end{equation}
On the one hand, applying Lemma \ref{recurrent} to the right-hand side of (\ref{eq4.7}), we  deduce that
\begin{equation*}
\nu\Big((f^{n-1})^{\ast}(f^{\ast})T,z\Big) =
\nu\Big((f^{n-n_0-1})^{\ast}[\mathcal{H}_0],z\Big)
 \end{equation*}
for any current $T:=[l(z)=0],$ where $\{l(z)=0\}$ is  a generic complex hyperplane in $\P^k,$ 
and for a generic point $z$ in any irreducible component of $ (f^{n-n_0-1})^{-1}(\mathcal{H}_0).$  On the other hand, under the same condition,
\begin{equation*}
\nu\Big((f^{n})^{\ast}T,z\Big) =0.
 \end{equation*}
We combine the latter two equalities with (\ref{eq4.5}) and taking into account the result of Step I.
Consequently, (\ref{eq4.6}) follows. This completes Step II.  The proof of the theorem is thereby finished.  \hfill $\square$

\begin{thm}\label{prop_main}
Let $f$ be a QAS meromorphic self-map of $\P^k.$ Then its first dynamical degree 
 $\lambda_1(f)$ is an algebraic integer.
Moreover,
\begin{equation}\label{eq_prop_0}
\lim\limits_{n\to\infty}\frac{ \d(f^{n+1})}{\d(f^n)}= \lambda_1(f).
\end{equation} 
\end{thm}
\begin{proof} 
If $f$ is AS, then  the theorem  is trivial since
$\d(f^n)=\d(f)^n$ and    $\lambda_1(f)=\d(f).$
Suppose now that  $f$ is non AS. Then 
in virtue of identities (\ref{eq4.1})--(\ref{eq4.2}), we have that
\begin{equation}\label{eq_prop_1}
\d(f^{n})=
 \begin{cases}
\d(f)^n, 
  &n=0,\ldots,n_0  ,\\
 \d(f)\cdot\d(f^{n-1})-\deg(H_0)\cdot \d(f^{n-n_0-1}), 
&  n>n_0.
\end{cases}
\end{equation}

Let  $h:=\deg(H_0).$ Then the characteristic polynomial
 associated to the sequence $ (\d(f^n))_{n=0}^{\infty}$ in  (\ref{eq_prop_1}) is
\begin{equation}\label{eq_prop_2}
P(\sigma):=\sigma^{n_0+1}-d\sigma^{n_0}+h.
\end{equation}
Let $\sigma_1,\ldots,\sigma_{s}$ be all the roots of $P$  with multiplicities
$m_1,\ldots,m_s$ respectively ($\sum\limits_{j=1}^s m_j=n_0+1$). 
Since  $P\in\R[\sigma],$ we remark that
  for every $p:\ 1\leq p\leq s,$ there is a unique $q$ such that
$\sigma_p=\overline{\sigma}_q$ and $m_p=m_q.$ 
Then the Newton formula yields
\begin{equation}\label{eq_prop_3}
\d(f^n)=\sum\limits_{j=1}^s P_j(n)\sigma^n_j,\qquad n\in\N,
\end{equation}
where $P_j$ is a complex polynomial of one variable  and $\deg(P_j)\leq m_j-1,$  $j=1,\ldots,s.$

There are two
cases  to consider.

\smallskip

\noindent{\bf Case 1:} all roots of $P$ are distinct.

In this case we know from (\ref{eq_prop_3}) that all polynomials $P_j$ appearing in
this formula are constant. Namely,  $P_j\equiv c_j,$  $j=1,\ldots,n_0+1.$ In virtue of
(\ref{eq_prop_1}) and (\ref{eq_prop_3}), we obtain the following
system of linear equations
with unknowns $c_1,\ldots,c_{n_0+1}:$
\begin{equation}\label{eq_prop_case1_1}
\left\lbrace
\begin{array}{l}
c_1+c_2+\cdots+c_{n_0+1}=1\\
\sigma_1c_1+\sigma_2 c_2+\cdots+\sigma_{n_0+1}c_{n_0+1}=d\\
\vdots  \\
\sigma_1^{n_0}c_1+\sigma_2^{n_0} c_2+\cdots+\sigma_{n_0+1}^{n_0}c_{n_0+1}=d^{n_0}\\
\end{array}
\right.
\end{equation}
Since all roots of $P$ are distinct, this system is nonsingular.

We  like to prove
the following\\
\noindent{\bf Claim 1.} {\it  If   $\sigma_q=\overline{\sigma_p}$  for some  $1\leq p<q\leq n_0+1,$ then
   $c_q=\overline{c_p}.$}

Indeed, using the remark made just after (\ref{eq_prop_2}), the assumption of Case 1 and solving system 
(\ref{eq_prop_case1_1}), Claim 1 follows.

\noindent{\bf Claim 2.} {\it  If  $\vert \sigma_p\vert=\vert \sigma_q\vert$  for some  $1\leq p<q\leq n_0+1,$ then
 $\sigma_q=\overline{\sigma_p}.$}  
 
  Indeed,
it follows from   (\ref{eq_prop_2}) and the identities  $P(\sigma_p)=P(\sigma_q)=0$ that
\begin{equation*}
h=\vert\sigma_p\vert^{n_0}\vert  d-\sigma_p\vert=\vert\sigma_q\vert^{n_0}\vert
d-\sigma_q\vert.
\end{equation*}
This, combined with the equality   $\vert \sigma_p\vert=\vert \sigma_q\vert,$ 
 implies that $\vert  d-\sigma_p\vert=\vert
d-\sigma_q\vert.$ In summary, we have that  $\vert  d-\sigma_p\vert=\vert
d-\sigma_q\vert$
and $\vert\sigma_p\vert=\vert\sigma_q\vert.$
Hence,    $\sigma_q=\overline{\sigma_p}.$ 
This completes the proof of  Claim 2.

We may assume without loss of generality that 
\begin{equation}\label{eq_prop_case1_2}
\vert
\sigma_1\vert=\max\limits_{c_j\not=0,\  j=1,\ldots, n_0+1}{\vert \sigma_j\vert}.
\end{equation} 
 Hence, $c_1\not=0.$
Next,  we show that   $\sigma_1$ is a real number. Suppose in order to get a contradiction that
this is not the case.  Without loss of generality we assume that   $\sigma_2=\overline{\sigma_1}.$
Then, by Claim 1, $c_2=\overline{c_1}.$
Write  $\sigma_1:=\vert\sigma_1\vert e^{i\phi}$ and $c_1=\vert c_1\vert e^{i\psi}$
for $0\leq \phi,\psi <2\pi,\ \phi>0$  and $\phi\not=\pi.$ 
Then formula   (\ref{eq_prop_3})  yields that
\begin{equation*}
\frac{\d(f^n)}{\vert \sigma_1\vert^n}=\frac{2\text{Re}(c_1\sigma_1^n)}{\vert\sigma_1\vert^n}
+\sum\limits_{j>2}c_j\left(\frac{\sigma_j}{
\vert\sigma_1\vert}\right)^n\\
=2\vert c_1\vert\cos{(n\phi+\psi)}+\sum\limits_{j>2}c_j\left(\frac{\sigma_j}{
\vert \sigma_1\vert}\right)^n,
\end{equation*}

Recall from  Claim 2  and (\ref{eq_prop_case1_2}) that  $\vert\sigma_1\vert  >\vert\sigma_j\vert$  for all $j>2$  with $c_j\not=0.$ 
Letting $n$ to infinity in the above formula, we see that any
accumulation point of the left side is a positive number, on the other
hand the sum in the right side tends to zero. Thus any accumulation point
of the set $\left\lbrace \cos{(n\phi+\psi)}\ :n\in\N\right\rbrace$
is also a positive number. We consider two subcases
\begin{itemize}
\item {\it (a)} $\frac{\phi}{2\pi}$ is rational.  Since  $0<\phi<2\pi$  and 
 $\phi\not=\pi,$
 there is a subsequence $( n_k)_{k= 1}^{\infty}$ such that $ \cos{(n_k\phi+\psi)}=c <0$
for all $k.$ This contradicts our observation above.

\item {\it (b)} $\frac{\phi}{2\pi}$ is irrational. It is clear that the
set $\left\lbrace\cos{(n\phi+\psi)}\ :n\in\N\right\rbrace$ is dense on
the interval $\lbrack-1,1\rbrack. $ Therefore  we obtain  a
contradiction.
\end{itemize}

In summary, we have shown that the algebraic integer $\sigma_1$ is a real number.
Applying Claim 2  we see that   $\vert \sigma_1\vert > \vert
\sigma_j\vert$ for any  $j\geq 2 $  with $c_j\not=0.$  Then it follows from (\ref{eq1.1}) and (\ref{eq_prop_3})   
 that $\sigma_1>0,$
$\lambda_1(f)=\sigma_1 $  and that  (\ref{eq_prop_0}) holds.
This completes the proof of Case 1.

 \noindent{\bf Case 2:} the polynomial $P$ has a multiple root.

Since  $h\geq 1,$ it can be checked that
$P(\sigma)=P^{'}(\sigma)=0$ if and only if $\sigma=\sigma_{n_0}:=
\frac{dn_0}{n_0+1}$  and  $h=h_0:=\left(\frac{d}{n_0+1}\right)^{n_0+1}n_0^{n_0}.$
Moreover, $P''(\sigma_{n_0})\not=0$ for $h=h_0.$  Consequently, we conclude that $h=h_0$ and the only multiple root of $P$ is 
$\sigma_{n_0},$ and it   is a double root.
In the remaining part of the proof we may assume  without loss of generality that $\sigma_{n_0+1}=\sigma_{n_0}.$

 Consider the system  of linear equations
with unknowns $c_1,\ldots,c_{n_0+1}:$
\begin{equation}\label{eq_prop_case2_1}
\left\lbrace
\begin{array}{l}
c_1+\cdots+c_{n_0-1}+c_{n_0+1}=1\\
\sigma_1c_1+\sigma_2 c_2+\cdots+\sigma_{n_0-1} c_{n_0-1}+ \sigma_{n_0}c_{n_0}+\sigma_{n_0}c_{n_0+1}=d\\
\vdots  \\
\sigma_1^{n_0}c_1+\sigma_2^{n_0} c_2+\cdots+\sigma_{n_0-1}^{n_0}c_{n_0-1} +        n_0\sigma_{n_0}^{n_0} c_{n_0}+      \sigma_{n_0}^{n_0}c_{n_0+1}=d^{n_0}\\
\end{array}
\right.
\end{equation}
and the formula
\begin{equation}\label{eq_prop_case2_2}
\d(f^n)=\sum\limits_{j=1}^{n_0-1} c_j\sigma^n_j+(nc_{n_0}+c_{n_0+1})\sigma_{n_0}^n,\qquad n\in\N.
\end{equation}

Using (\ref{eq_prop_case2_1})--(\ref{eq_prop_case2_2})
instead of (\ref{eq_prop_3}) and  (\ref{eq_prop_case1_1}),
and taking into account that  $\sigma_{n_0+1}=\sigma_{n_0}$ is a positive number, we can proceed     as in Case 1.
There is only a small caution. Namely, the corresponding  Claim 1 should be read as follows:\\
 {\it  If   $\sigma_q=\overline{\sigma_p}$  for some  $1\leq p<q\leq n_0-1,$ then
  $c_q=\overline{c_p}.$}

Let 
\begin{equation} \label{eq_prop_case2_3} 
\lambda:=\max\limits_{c_j\not=0,\  j=1,\ldots, n_0+1}{\vert \sigma_j\vert}.
\end{equation} 
There are two subcases to consider.\\
\noindent {\bf Subcase 2a:} $\lambda\not=\sigma_{n_0}.$

Then  $n_0>1.$ We may suppose without loss of generality that  $\vert \sigma_1\vert =\lambda.$  Using  (\ref{eq_prop_case2_2})--(\ref{eq_prop_case2_3}),   
     we proceed as in Case 1. Hence, the desired conclusion follows.\\
\noindent {\bf Subcase 2b:} $\lambda=\sigma_{n_0}.$
 
Using  (\ref{eq_prop_case2_2})--(\ref{eq_prop_case2_3}),  we see that 
\begin{equation*}
\frac{\d(f^n)}{\vert \sigma_{n_0}\vert^n}=(nc_{n_0}+c_{n_0+1})+
\sum\limits_{j<n_0}c_j\left(\frac{\sigma_j}{\vert\sigma_{n_0}\vert}\right)^n.
\end{equation*}
Since  the sum in the right side tends to zero as $n\to\infty,$  it follows  that  $\lambda_1(f)=\sigma_{n_0}$
and that  (\ref{eq_prop_0}) holds. Hence, the proof of Case 2 follows.

 This finishes the proof. 
\end{proof}

 \subsection*{Applications.}  We  apply  the Main Theorem  to the examples given in Section 3.
 
 First   consider Example  \ref{ex1}.
  Put $H(z,w,t):=t$ and  let $F:\ \C^3\longrightarrow\C^3$  be  the lifting of $f$ given  by the right-hand
side of (\ref{eqex1}).
 In the light of the Main Theorem,
 a sequence of liftings $F_n$ of $f^n$  may be defined  as follows
\begin{equation*}
F_0:=\Id,\ F_1:= F,\qquad F_{n}:=\frac{F\circ F_{n-1}}
{     H(F_{n-2})
},\qquad  n\geq 2.
\end{equation*}
  As a consequence, one obtains  the  equation   $\d(f^n)-2\d(f^{n-1})+\d(f^{n-2})=0.$   Therefore,
 a straightforward computation shows that
$\d(f^n)=n+1 $ and  $\lambda_1(f)=1.$

 Next consider Example  \ref{ex3}.
  Put $H(z,w,t):=(z+w+t)^2(z^3+w^3+t^3)$ and  let $F:\ \C^3\longrightarrow\C^3$  be  the lifting of $f$ given by  the right-hand
side of (\ref{eqex3}).
 Using  the Main Theorem,
 a sequence of liftings $F_n$ of $f^n$  may be defined  as follows
\begin{equation*}
F_0:=\Id,\ F_1:= F,\qquad F_{n}:=\frac{F\circ F_{n-1}}
{     H(F_{n-2})
},\qquad  n\geq 2.
\end{equation*}
  As a consequence, one obtains  the  equation   $\d(f^n)-7\d(f^{n-1})+5\d(f^{n-2})=0.$   Therefore,
 a straightforward computation shows that
\begin{equation*}
\d(f^n)=\frac{ 1}{\sqrt{29}}\cdot\Big(\frac{7+\sqrt{29}}{2}\Big)^{n+1}- \frac{ 1}{\sqrt{29}}\cdot\Big(\frac{7-\sqrt{29}}{2}\Big)^{n+1} ,
\end{equation*} and  $\lambda_1(f)=\frac{7+\sqrt{29}}{2}.$
\subsection*{Concluding remarks.}  One  may widen  the class of QAS self-maps by weakening  considerably the conditions in Definition \ref{QAS}.
Of course the recurrent law would be then more complicating.
  One might also seek to 
  \begin{itemize}
  \item[$\bullet$] for any given $k,d\geq 2,$ find   many
families of QAS (but non AS) self-maps  of $\P^k$ with the algebraic degree $d;$
 \item[$\bullet$] generalize  the Main Theorem  to meromorphic self-maps in  compact K\"{a}hler manifolds;
  \item[$\bullet$] construct an {\it appropriate} Green current for  every QAS self-map $f$ with $\lambda_1(f)>1.$
  \end{itemize}
  
  We hope to come back  these issues in a future work.

\end{document}